\documentclass[11pt,a4paper]{article}
\usepackage[utf8]{inputenc}
\usepackage[english]{babel}
\usepackage[centertags]{amsmath}
\usepackage{amsfonts,amssymb,amsthm}
\usepackage{xcolor}
\usepackage{hyperref}
\usepackage{multicol}
\newtheorem{theorem}{Theorem}
\newtheorem{lemma}{Lemma}

\newtheorem{proposition}{Proposition}
\newtheorem{remark}{Remark}
\newtheorem{corollary}{Corollary}
\newcommand{\N}{\mathbb N}
\newcommand{\R}{\mathbb R}
\renewcommand{\L}{\mathcal L}
\newcommand{\K}{\mathcal K}
\newcommand{\cN}{\mathcal N}
\newcommand{\cNi}{\tilde{\mathcal N}}
\newcommand{\eps}{\varepsilon}
\newcommand{\E}{\mathbb{E}\;}
\newcommand{\ind}{\mathbb{I}\;}
\newcommand{\var}{{\rm{var}}\;}
\newcommand{\cs}{{\rm{cst}}\;}
\renewcommand{\P}{\mathbb P}
\newcommand{\<}{\left\langle}
\renewcommand{\>}{\right\rangle}

\title{On the number of roots of Sturm-Liouville random sums}
\date{\today}
\author{Federico Dalmao\thanks{DMEL, Cenur Litoral Norte, Universidad de la República, Uruguay,  fdalmao@unorte.edu.uy.}\hspace{2mm}  
 \& José R. León\thanks{IMERL, FING, Universidad de la República, Uruguay 
 and Escuela de Matemática, Universidad Central de Venezuela, rlramos@fing.edu.uy.}}
\begin{document}
\maketitle
\begin{abstract}
We consider the number of roots of linear combinations of a system of $n$ orthogonal eigenfunctions of a Sturm-Liouville 
initial value problem with i.i.d. standard Gaussian coefficients. 
We prove that its distribution inherits the asymptotic behavior of the number of roots 
of Quall's random trigonometric polynomials. 
This result can be thought as a robustness result for the central limit theorem for the number of roots of Quall's random trigonometric polynomials 
in the sense that small uniform perturbations of sines and cosines do not change the limit distribution.
\end{abstract}

\section{Introduction} 
The study of the number of zeros of random polynomials 
is an active field of research in mathematics and physics. 
Within mathematics it lies at the intersection 
of analysis, probability, algebraic geometry 
among other areas. 
The reader can profit of a centenary literature with 
a constellation of different algebraic structures, 
kinds of randomness sources and techniques, 
see e.g. \cite{bharucha} and \cite{Farahmand} and references therein. 

One of the main classes of random polynomials are 
\emph{random trigonometric polynomials} of the forms : 
\begin{equation} \label{d:RTP}
 T_n(x) := \frac{1}{\sqrt{n}}\sum^n_{k=1}a_k\cos(kx)+b_k\sin(kx) 
 \qquad\rm{or}\qquad 
 C_n(x) := \frac{1}{\sqrt{n}}\sum^n_{k=1}a_k\cos(kx), 
\end{equation}
where the coefficients are standardized r.v.s and $t$ is restricted to some subinterval of $[0,2\pi]$.  
The literature concerned with their zeros have been increasing quickly 
in the past years, see e.g. \cite{adl,al,pautrel,vietnam1,vietnam} and references therein. 
We refer to $T_n$ as Quall's or stationary random trigonometric polynomials 
and to $C_n$ as classical random trigonometric polynomials.

A trigonometric polynomial can be thought as a finite Fourier expansion 
of a given function. 
Thus, random trigonometric polynomials are, in this sense, finite Fourier expansions of random functions, that is, 
expansions in the orthogonal basis formed by sines and cosines. 
Naturally, one can think of other similar expansions as that in terms of Bessel functions 
or in terms of other ensembles of orthogonal polynomials on an interval, see e.g. \cite{Lubinsky,vietnam}. 
Fourier and Bessel expansions are particular cases of 
Sturm-Liouville expansions. 
To introduce them, 
let $q:[a,b]\to\R$ be a positive or bounded by below and continuous function 
with finite limits at $a,b$ 
and consider the differential operator
\begin{equation} \label{d:SLdop}
 \L := q(x) - \frac{d^2}{dx^2},
\end{equation}
which acts on smooth functions defined on $[a,b]$. 
Under these conditions on $q$, the problem is called regular. 
When $q=0$ we retrieve the trigonometric/Fourier case \cite[pp8,pp22]{Titchmarsh} while if $q(x)=(\nu-\frac14)\frac{1}{x^2}$ we recover the Bessel case 
\cite[S.1.11]{Titchmarsh}.

The classical theory, states that 
under mild conditions, see \cite[Th.1.9]{Titchmarsh}, 
an integrable function $f$ can be expanded in terms of the eigenfunctions 
of $\L$, 
being the behavior of this expansion similar to that of the Fourier case. 
A function $\psi$ is said to be an eigenfunction of $\L$ if 
\begin{equation} \label{e:eigen}
	\L\psi=\lambda\psi,
\end{equation}
for some number $\lambda$ (known as the corresponding eigenvalue) and provided that it verifies some specified border conditions. 
See the details below. 
In particular, it is well known, see \cite[S1.12]{Titchmarsh}, 
that the eigenvalues are simple, positive and unbounded. That is, the eigenvalues form a sequence
\[
0<\lambda_1<\cdots<\lambda_n <\cdots \to \infty.
\]
Moreover, the eigenfunctions $\psi_n$ of $\L$, corresponding to the eigenvalues $\lambda_n:n\in\N$, form an orthogonal system, thus generalizing the trigonometric case.

\begin{remark}[Generalizing the Wave Equation]
One of the paradigmatic applications of Fourier expansions 
is the solution of the (here, one dimensional) \emph{Wave Equation}:
\[
 \dfrac{\partial^2 u}{\partial t^2} 
= c\dfrac{\partial^2 u}{\partial x^2}, \quad t\in[a,b],
\] 
for some constant $c>0$. 
The eigenfunctions of a differential operator related to $\L$ help to solve  more general initial value problems of the forms
\begin{equation}\label{e:de}
(C): 
\begin{cases}
\dfrac{\partial^2 u}{\partial t^2} 
= \dfrac1{\omega^2(x)}\dfrac{\partial^2 u}{\partial x^2}\\
u(t,a)=u_a(t), \ t>0,\\
 u'(t,a)=0,\ t>0,
\end{cases}
\qquad\qquad
(D) : 
\begin{cases}
\dfrac{\partial^2 u}{\partial t^2} 
= \dfrac1{\omega^2(x)}\dfrac{\partial^2 u}{\partial x^2}\\
u(t,a) = u(t,b)=0,\ t>0,\\
\end{cases}
\end{equation}
where $\omega$ is a strictly positive, twice continuously differentiable function s.t. $\int^b_a w(u)du<\infty$.

Indeed, in case $(D)$, let $\K:f\mapsto \frac1{\omega^2(x)} f''$, 
we say that a non trivial function $f$ is an eigenfunction of $\K$ with eigenvalue $\gamma$ if $\K (f) = \gamma f$ and $f(a)=f(b)=0$.  
Thus, we deduce that 
$f''(x)=\gamma \omega^2(x)f(x)$, for $x\in[a,b]$. 
If we consider the Hilbert space $\mathbb L^2 \big([a,b],\omega^2(x)dx \big)$ and denoting $\<,\>_{\omega^2}$ 
its scalar product, we have
\[
\<\K f,g\>_{\omega^2} =\int_a^b\frac1{\omega^2(x)}f''(x)g(x)\omega^2(x)dx =\<f,\K g\>_{\omega^2}.
\] 
In this form, the operator is self-adjoint and moreover 
if $f_\gamma$ is an eigenfunction associated to the  $\gamma$ eigenvalue we have
\[
\<\K f_\gamma,f_\gamma\>_{\omega^2} 
=\gamma||f_{\gamma}||^2_{\omega^2}
\quad {\rm and }\ 
\<\K f_\gamma,f_\gamma\>_{\omega^2} 
=\int_a^b f''_\gamma(x)f_\gamma(x)dx
=-\Big\|\frac{f'_\gamma}{\omega}\Big\|^2_{\omega^2}.
\]
This implies that $\gamma<0$, if it exists. 
The operator $\K$ is self-adjoint and closed. 
We claim that the subspace ${\rm Ker}(\K)=\{0\}$. 
In fact, if $f\in{\rm Ker}(\K)$ it holds $\frac1{\omega^2(x)}f''(x)=0$, 
then $f(x)=c_1+c_2 x$, $x\in[a,b]$, 
but the boundary condition implies that $c_1=c_2=0$ and the claim follows. 
In particular, all the eigenvalues are different from zero. 
The general theory \cite[Ch.7]{Coddington} and \cite{BR} gives that the eigenspaces are of dimension one and that the 
spectrum is discrete $\{\gamma_n\}_{n=1}^\infty$.  
If we denote $f_n$ the eigenfunctions of $\K$ we obtain
that the general solution to case $(D)$ in \eqref{e:de} can be written as
\[
u(t,x)=\sum_{n=1}^\infty \Bigg[a_n\cos\big((-\gamma_n)^{\frac12}t\big) + b_n\frac{ \sin\big( (-\gamma_n)^{\frac12}t\big)} {(-\gamma_n)^{\frac12}}\Bigg] f_n(x).
\]
Case $(C)$ and more general cases can be treated similarly 
see also \cite{Coddington}.\hfill$\triangleleft$
\end{remark}

\begin{remark}[Normal form of a Sturm Liouville system] \label{r:normal}
The operator $\L$ in \eqref{d:SLdop} corresponds to the so called \emph{normal form} and the operator $\K$ associated to the initial value problems \eqref{e:de} can be translated to the normal form as follows. 

From now on, we put $a=0$ and $b=2\pi$ 
and we assume w.l.o.g. that $\int^{2\pi}_0\omega(u)du=2\pi$. 
First observe that an eigenfunction $\psi$ of $\L$ verifies 
\begin{equation*} %\label{titchmersh:141}
\psi''-q\psi=\lambda \psi,
\end{equation*}
see \cite[Eq.(1.1.4)]{Titchmarsh}, \cite[Eq.s(2.2)-(2.4)]{Fulton} or
\cite[Eq.s(39)-(40)]{BR}. 
Now, take an eigenfunction $f$ of $\K$, consider the following change of variables 
$y= \int_0^x\omega(u)du$ 
and define $g(y) :=  \omega^{1/2}(x)f(x)$, $x\in[0,2\pi]$. 
We have $\frac{dx}{dy}=\frac{1}{\omega(x)}$. 
Thus, omitting $(x)$, we have
\[
\frac{dg}{dy} = \frac{dg}{dx}\frac{dx}{dy}
	=\frac{dg}{dx}\frac{1}{\omega} 
	=\frac{1}{\omega} \Big[\frac{\omega' f}{2\omega^{1/2}} 
	+ \omega^{1/2} f' \Big] 
	= \frac{\omega' f}{2 \omega^{3/2}} + \frac{f'}{\omega^{1/2}}.
\]
Then,
\begin{equation*}
\frac{d^2g}{d^2y} = \frac{1}{\omega} 
\frac{d}{dx} \Big[ \frac{\omega' f}{2 \omega^{3/2}} 
+ \frac{f'}{\omega^{1/2}} \Big] 
= \Big[ \frac{\omega''}{2\omega^{3}} 
- \frac34\frac{(\omega')^2}{\omega^{4}} 
 + \gamma  \Big]\cdot g.
\end{equation*}
In the last equality we used that $f''=\gamma\omega^2 f$ 
and that $g = \omega^{1/2}f$. 
Thus, the resulting equation 
$\frac{d^2g}{d^2y} = \Big[\gamma + \frac{\omega''}{2\omega^3} 
-\frac{3(\omega')^2}{4\omega^{4}}\Big] g$ 
matches \eqref{e:eigen}, see also \cite[Eq.(1.1.4)]{Titchmarsh}, 
with $q=\frac{\omega''}{2\omega^{3}} 
- \frac34\frac{(\omega')^2}{\omega^{4}}$ and $\lambda=-\gamma$.  
\hfill$\triangleleft$
\end{remark}
 
%%%%%
The object of this work is the study of the number of zeros of 
random (Gaussian) Sturm-Liouville sums, that is, random linear combinations 
of eigenfunctions of $\L$ of the form.
\begin{equation*}
	\frac{1}{\sqrt{n}} \sum^n_{k=1} a_k\psi_k(x),\ x\in[a,b],
\end{equation*}
with i.i.d. standard Gaussian $a_k:k\geq 1$, 
as $n\to\infty$. 
See \eqref{d:Fn} for the precise definition. 
It is a well known fact
that the eigenfunction of regular Sturm-Liouville problems 
can be approximated by trigonometric functions. 
This suggests that the number of zeros, or even the zeros themselves, 
of random Sturm-Liouville sums should be close to those 
of random trigonometric polynomials. 
Making this idea precise is one of the goals of the present work. 
One is tempted to repeat the arguments through chaotic decompositions 
which led to central limit theorems for the number of zeros of random trigonometric polynomials, 
see \cite{adl,al} and also \cite{vietnam}, 
but, on the other hand the additional information should allow, 
may be through Rice's formulas, 
to inherit the asymptotic distribution for the number of zeros of random Sturm-Liouville sums 
from that of the number of roots of random trigonometric polynomials. 
In the stronger a.s. sense, 
the recent closed (i.e.: without the limit) Kac-type formula obtained in \cite[Pr.1]{ap} is appealing. 

We prove, using classical Kac' formula, 
that the asymptotic behavior, as $n\to\infty$, 
of the number of zeros of random Sturm-Liouville sums 
coincide with that of random trigonometric polynomials, 
that is, its asymptotic mean and variance are of order $n$ and, after standardization, its distribution converges towards a centered normal one. 

Let us finish the introduction 
pointing out that our result can be viewed as a robustness result 
since it implies that if one perturbs, see the details below, sines and cosines 
one still has the same limit distribution for the number of zeros. 
\medskip

The paper is organized as follows. 
We set the problem and state our main result in Section 2. 
In Section 3 we consider some important preliminary facts. 
The proof of the main result is presented in Section 4. 
Some anciliary computations are deferred to Section 5.

Regarding the notation, $\cs$ stands for a positive constant whose value may change at each appearance; 
$a_n=O(b_n)$ means that $|a_n|\leq \cs b_n$ for $n$ large enough;
$a_n\sim_nb_n$ means that $\lim_{n\to\infty}\tfrac{a_n}{b_n}=1$.

\section{Problem setting and main result} 
Let $\omega$ be a strictly positive, twice continuously differentiable weight function s.t. $\int^{2\pi}_0\omega(u)du=2\pi$. 
Consider a sequence of independent standard Gaussian 
r.v.s $\{a_k, b_k\}_{k=1}^\infty$ 
and define 
\begin{equation} \label{d:Fn}
F_n(x) := \frac{1}{\sqrt{n}} \sum_{k=1}^n a_k u_k(x) + b_k v_k(x),\ x\in[0,2\pi],
\end{equation}
where $u_k$, $v_k : k\geq 1$ stand for sequences 
of orthonormal eigenfunctions (with eigenvalue $\lambda_k$)
of the Sturm-Liouville operator $\L$, \eqref{d:SLdop}-\eqref{e:eigen}, associated to $q=\frac{\omega''}{2\omega^{3}} 
- \frac34\frac{(\omega')^2}{\omega^{4}}$ (see Remark \ref{r:normal}), corresponding respectively to 
two sets of 'basic'{} initial conditions 
\begin{align*}
(C)   :\ g(0) \neq 0, g'(0) = 0,\qquad
(D)  :\ g(0)  = g(2\pi) = 0,
\end{align*}
whose existence is guaranteed by \cite[Th.1.5]{Titchmarsh}.

Let $\cN_n$ be the number of roots of $F_n$ on the interval $[0,2\pi]$, 
i.e:
\begin{equation*} %\label{d:Nn}
 \cN_n := \#\{x\in[0,2\pi] : F_n(x)=0\}.
\end{equation*} 
The next theorem is the main result of this work.
\begin{theorem} \label{t:ppal}
With the above notation, as $n\to\infty$, there exists $0<V<\infty$ s.t.
\begin{equation*}
 \lim_{n\to\infty} \frac{\var(\cN_n)}{n} =V,
\end{equation*}
and after normalization, the distribution of $\cN_n$ converges towards the standard Gaussian law.
\end{theorem}

The main idea of the proof of Theorem \ref{t:ppal} 
is to take profit of the available central limit theorem for the number of zeros of 
stationary trigonometric polynomials in \cite{gw} and \cite{al} by 
assesing the $L^1$ contiguity between both numbers of zeros.
As a by-product of our proof 
we obtain the following \emph{robustness} result for \emph{perturbed} random trigonometric polynomials. 
\begin{corollary}
 Let $\eps_k,\eta_k: [0,2\pi]\to\R$, $k\geq 1$ be of class $C^2$ with $|\eps_k(\cdot)|, |\eta_k(\cdot)|\leq\tfrac{\cs}{k}$ 
 and $|\eps'_k(\cdot)|, |\eta'_k(\cdot)|\leq\cs$\hspace{-1.4mm}, 
 then, the conclusions of Theorem \ref{t:ppal} hold true for the number of roots of the 
 perturbed random trigonometric polynomial
 \[
    \frac{1}{\sqrt n}\sum^n_{k=1} a_k (\cos(kx)+\eps_k(x)) +b_k (\sin(kx)+\eta_k(x)),
 \] 
 for i.i.d standard Gaussian $a_k,b_k:k\geq 1$. 
\end{corollary}

\begin{remark}
 Here, we target stationary random trigonometric polynomials, 
 but, it is clear that one can consider random sums involving $\{u_K:k\geq 1\}$ only for instance, 
 which shall be approximated by classical trigonometric polynomials \cite{adl}.
\end{remark}

\section{Preliminaries}
It is well known that 
the eigenvalues of the operator $\L$ are real 
and that they verify the asymptotics 
$\sqrt{\lambda_n} = \frac{n}{2}+O(\frac{1}{n})$ in case (C) 
and $\sqrt{\lambda_n} = \frac{n+1}{2}+O(\frac{1}{n})$ in case (D), 
see \cite[S.1.12]{Titchmarsh} or \cite{BR} and \cite{Fulton}. 
For ease of notation we group the eigenfunctions according to their phase $\frac{n}{2}$.

Besides, Lemmas 1.7 (ii) and (iii) of \cite{Titchmarsh} states that  
the eigenfunctions of $\L$ verify the asymptotics 
(w.l.o.g. we choose the leading constant to be $1$, 
afterwards we will normalize the polynomials so that 
they have variance one):
\begin{itemize}
\item in case (C) : for $n$ large enough, 
uniformly in $[0,2\pi]$, we have
\begin{align*} %\label{e:aprox-un}
u_n(x) &= \frac{1}{\omega(x)^{1/2}} \cos\left(\frac{n}{2}\int^x_0\omega(u)du\right) 
+ O\Big(\frac{1}{n}\Big),\\
u'_n(x) &= -\frac{n}{2} \omega(x)^{1/2} \sin\left(\frac{n}{2}\int^x_0\omega(u)du\right) 
+ O(1). \notag
\end{align*}

\item In case $(D)$, for $n$ large enough, 
uniformly in $[0,2\pi]$, we have
\begin{align*} %\label{e:aprox-vn}
v_n(x) &= \frac{1}{\omega(x)^{1/2}} \sin\left(\frac{n}{2}\int^x_0\omega(u)du\right) 
+ O\Big(\frac{1}{n}\Big), \\
v'_n(x) &= \frac{n}{2} \omega(x)^{1/2} \cos\left( \frac{n}{2}\int^x_0\omega(u)du\right)
+ O(1). \notag
\end{align*}
\end{itemize}

\section{Proof of Theorem \ref{t:ppal}}
For a process $Z_n$ defined on the interval $I$, 
denote its number of zeros 
and its standardized number of zeros 
respectively by
\begin{equation} 
 \cN (Z_n,I) := \#\big\{ x\in I : Z_n(x) = 0\big\} \quad {\rm and }\quad
 \cNi (Z_n,I) := \frac{\cN(Z_n,I) - \E \cN(Z_n,I)}{\sqrt{n}}.\label{d:Nn}
\end{equation}

For $x\in[0,2\pi]$, 
denote $\Omega(x) := \int^x_0\omega(u)du$,
and 
\begin{equation} \label{e:c-s}
 c_k := \frac{1}{\omega(x)^{1/2}} \cos\left(\frac{k}{2} \Omega(x)\right), \qquad
 s_k := \frac{1}{\omega(x)^{1/2}} \sin\left(\frac{k}{2} \Omega(x)\right).
\end{equation}
Let also 
\begin{equation*} %\label{e:X}
X^o_{n}(x) 
=\frac{1}{\sqrt n} \sum_{k=1}^n a_k c_k(x) + b_k s_k(x),
\quad x\in[0,2\pi].
\end{equation*}

The next lemma gives the first approximation 
between the numbers of zeros of two processes, 
relating $X^o_n$ with the stationary trigonometric polynomials 
$T_n$ in \eqref{d:RTP} 
studied in \cite{gw} and \cite{al}. 
\begin{lemma}
The distribution of the (standardized) number of zeros $\cNi(X^o_n,[0,2\pi])$ of the process $X^o_n$ converges, as $n\to\infty$, towards a centered normal r.v. with finite and positive variance.  
\end{lemma}
%%% 
\begin{proof}
Since $\Omega:[0,2\pi]\to[0,2\pi]$ is bijective, 
setting $y=\Omega(x)$ and introducing the process
$Y^o_n(y)= X^o_n(\Omega^{-1}(y))$, we have
\begin{equation*}
Y^o_n(y)=\sqrt{ \frac{ 2}{n \omega(\Omega^{-1}(y))}} \sum_{k=1}^n a_k \cos\Big(\frac{k}{2} y\Big) + b_k \sin\Big(\frac{k}{2} y\Big).
\end{equation*}
As we are looking at the roots of these processes, 
the normalization $\sqrt{2/\omega(\Omega^{-1}(y))}$ plays no role. 
Besides, for $T\in[0,2\pi]$,
\[
\cN (X^o_n,[0,T]) := \{ x\in[0,T]:\, X^o_n (x) = 0\}
= \{y\in[0,\Omega(T)]:\, Y^o_n(y)=0\} 
=: \cN(Y^o_n,[0,\Omega(T)]).
\]
In particular, as $\Omega(2\pi)=2\pi$, 
we have $\cN(X^o_n,[0,2\pi]) = \cN(Y^o_n,[0,2\pi])$. 
Thus, we are led to study the number of roots of a \emph{stationary trigonometric polynomial} $T_n$ (see \eqref{d:RTP}) restricted to the interval 
$[0,\pi]$. 
The result follows.
\end{proof}

The next step is to approximate $\cN_n=\cN(F_n,[0,2\pi])$ by $\cN(X^o_n,[0,2\pi])$. 
It is convenient to standardize the processes $X^o_n$ and $F_n$. 
Define for $x\in[0,2\pi]$
\begin{equation} \label{d:Xn}
X_n(x) := \omega(x)^{1/2} X^o_n(x) 
= \frac{1}{\sqrt n} \sum_{k=1}^n 
a_k \cos\Big(\frac{k}{2}\Omega(x)\Big) 
+ b_k \sin\Big(\frac{k}{2}\Omega(x)\Big),
\end{equation}
and
\begin{equation} \label{d:fn}
 f_n(x) := \sqrt{\omega(x)} F_n(x).
\end{equation}
Observe that, since the factor $\sqrt{\omega(x)}$ in the definition of $f_n$ in \eqref{d:fn} plays no role in the study of the zeros, 
we have $\cN_n = \cN(f_n,[0,2\pi])$. 

The next proposition provides the final approximation we need. 
Once its established, the central limit theorem for $\cN(f_n,[0,2\pi])$ 
follows from that for $\cN(X_n,[0,2\pi])$. 
\begin{proposition} \label{p:NL2}
For $X_n$, $f_n$ defined as in \eqref{d:Xn} and \eqref{d:fn} 
and for $\cN(f_n,[0,2\pi])$, $\cN(X_n,[0,2\pi])$ defined as in \eqref{d:Nn}, 
we have 
\begin{equation*}
 \frac{\cN(f_n,[0,2\pi])-\cN(X_n,[0,2\pi])} {\sqrt n}\ 
 \mathop{\longrightarrow}\limits_{n}^{L^1}\ 0.
\end{equation*}
\end{proposition}
\begin{proof}
For brevity, we set $\cN_{f_n}$ and $\cN_{X_n}$ for $\cN(f_n,[0,2\pi])$ and $\cN(X_n,[0,2\pi])$ respectively.

We use the Kac formula to estimate the $L^1$ distance. 
We have
\begin{align*}
 \E|\cN_{f_n}-\cN_{X_n}| &= 
 \E\Big|\lim_{\delta\downarrow0}\frac{1}{2\delta} \int^{2\pi}_0 
 \big[ |X'_n| \ind_{|X_n|<\delta} - |f'_n| \ind_{|f_n|<\delta}\big]dx \Big| \\
% &= \E\lim_{\delta\downarrow0}\frac{1}{2\delta} \Big|\int^{2\pi}_0 
% \big[ |X'_n| \ind_{|X_n|<\delta} - |f'_n| \ind_{|f_n|<\delta} 
% \pm |X'_n| \ind_{|f_n|<\delta}\big]dx \Big| \\
 &\leq \lim_{\delta\downarrow0}\frac{1}{2\delta} \E \int^{2\pi}_0 
 \big[ |X'_n-f'_n| \ind_{|f_n|<\delta} + |X'_n| \big|\ind_{|X_n|<\delta}-\ind_{|f_n|<\delta}\big|\big]dx \\
 &=: A(n) + B(n).
\end{align*}
We used Fatou's and triangular inequalities to get the bound.

First, consider $A(n)$.  
Set $\alpha_n(x)=\frac{\E(X'_n(x)f_n(x))}{\var(f_n(x))}$, 
thus
$\tilde{X}'_n(x) := X'_n(x) -\alpha_n(x)f_n(x)$ is independent from $f_n(x)$. 
Besides, $|\alpha_n(x)|\leq\cs$ for 
$x\in[0,2\pi]$, see part 1 in Lemma \ref{l:aux}.
Then, we have
\begin{multline*}
A(n)= \lim_{\delta\downarrow0}\frac{1}{2\delta} \E \int^{2\pi}_0 
  |X'_n(x)-f'_n(x)| \ind_{|f_n(x)|<\delta}\,dx\\
  = \lim_{\delta\downarrow0}\frac{1}{2\delta} 	\E \int^{2\pi}_0 
  |\tilde{X}'_n(x)+\alpha_nf_n(x)-f'_n(x)| \ind_{|f_n|<\delta}\,dx
  =\E|X'_n(0)-f'_n(0)|\phi(0).
\end{multline*}
But, the random variable $X'_n(0)-f'_n(0)$ is Gaussian and furthermore
\[
 X'_n(0)-f'_n(0)=\frac{\omega(x)}{2\sqrt n}\sum_{k=1}^n 
 k\Big[a_kO\Big(\frac1k\Big)+b_kO\Big(\frac1k\Big)\Big]
 +\frac{\omega'(x)}{2\omega(x)}f_n(x).
\] 
The first term in the above sum is $N(0,O(1))$ and the second term is $N(0,1)$. 
This notation means that the first term has a bounded variance. These two facts imply $\E|X'_n(0)-f'_n(0)|=O(1)$, 
concluding that $\frac{A(n)}{\sqrt n}\to0$.

We move to $B(n)$. 
Hence, with the same notations as above, 
we have
\begin{align*}
 B(n) &= \lim_{\delta\downarrow0}\frac{1}{2\delta} \int^{2\pi}_0 \E \Big[
 |\tilde{X}'_n+\alpha_nf_n| \Big|\ind_{|X_n|<\delta}-\ind_{|f_n|<\delta}\Big|\big] \Big] dx \\
 &\leq \lim_{\delta\downarrow0}\frac{1}{2\delta} \int^{2\pi}_0 
 \big[ \E |\tilde{X}'_n|+\alpha \E|f_n|\big] 
 \P\Big\{\{ |X_n|<\delta, |f_n|>\delta\}\cup\{|X_n|>\delta, |f_n|<\delta\} \Big\} dx.
\end{align*}
We have $\E |\tilde{X}'_n| \sim_n \E |X'_n| = \cs n$ and $\E |f_n| \sim_n \cs$\hspace{-1mm}, see Lemma \ref{l:cov}. 
We consider one of the two terms, the other one is analogous. 
Denote 
\begin{equation*} 
  \eps_n(x) := f_n(x)-X_n(x). 
\end{equation*}
We begin with the help of a control in $\|\eps_n\|_\infty:=\sup_{x\in[0,2\pi]}|\eps_n(x)|$. 
We fix $x$, we omit the $(x)$ in $X_n(x)$, $f_n(x)$ and $\eps_n(x)$.
\begin{align*}
& \frac{1}{2\delta} \P\Big\{|X_n|<\delta, |f_n|>\delta, \|\eps_n\|_{\infty} <\frac{\cs\log n}{\sqrt n} \Big\} 
 \leq \frac{1}{2\delta} \P\Big\{|X_n|<\delta < |f_n| < \delta + \frac{\cs\log n}{\sqrt n} \Big\} \\
 &\qquad= \frac{1}{2\delta}\int^\delta_{-\delta} du \int_\delta^{\delta + \frac{\cs\log n}{\sqrt n}}  p_{{}_{X_n,f_n}} (u,v) dv 
 \ \mathop{\to}\limits_{\delta\downarrow0}\ \int_0^{\frac{\cs\log n}{\sqrt n}}  p_{{}_{X_n,f_n}} (0,v) dv \\
 &\qquad =\ \frac{\cs\log n}{\sqrt n} p_{{}_{X_n,f_n}} \Big(0,\theta \frac{\cs\log n}{\sqrt n}\Big) 
 = \frac{\cs\log n}{\sqrt n} \frac{1}{2\pi\sqrt{\Delta}} \exp\Big\{-\frac12 \frac{\var(f_n)\cdot (\frac{\cs\log n}{\sqrt n})^2}{\Delta}\Big\}.
\end{align*}
Here $\Delta$ stands for the determinant of $\var(X_n(x),f_n(x))$.
As $\Delta \sim_n \frac{\cs}{n}$, see part 2 in Lemma \ref{l:aux},
we have
\begin{equation*}
 \lim_{\delta\downarrow0} \frac{1}{2\delta} 
 \P\Big\{|X_n|<\delta, |f_n|>\delta, \|\eps_n\|_{\infty} <\frac{\cs\log n}{\sqrt n} \Big\} 
 \leq \cs\log n\cdot \exp\{-\cs(\log n)^2\} = \frac{\cs\log n}{n^{\cs\log n}}.
\end{equation*}

Now, we deal with the remainder. 
Fix $x\in[0,2\pi]$ for a moment.
\begin{multline*}
 \frac{1}{2\delta} \P\Big\{|X_n(x)|<\delta, |f_n(x)|>\delta, \|\eps_n\|_{\infty} \geq \frac{\cs\log n}{\sqrt n} \Big\} 
 \leq \frac{1}{2\delta} \P\Big\{|X_n(x)|<\delta, \|\eps_n\|_{\infty} \geq \frac{\cs\log n}{\sqrt n} \Big\} \\
 = \frac{1}{2\delta}\int^\delta_{-\delta} \P\Big\{ \|\eps_n\|_{\infty}\geq \frac{c\log n}{\sqrt n}\mid X_n(x) = u \Big\} 
 p_{{}_{X_n(x)}}(u)du.
\end{multline*}
For $y\in[0,2\pi]$, set $\beta_n(y)$ s.t. $\eps_n(y)-\beta_n(y)X_n(x)$ is independent from $X_n(x)$. 
Thus, 
\begin{align*}
 &\frac{1}{2\delta} \P\Big\{|X_n(x)|<\delta, |f_n(x)|>\delta, \|\eps_n\|_{\infty} \geq \frac{\cs\log n}{\sqrt n} \Big\} \\
 &\qquad = \frac{1}{2\delta}\int^\delta_{-\delta} \P\Big\{ \|\eps_n(\cdot)-\beta_n(\cdot)X_n(x) +\beta_n(\cdot) u\|_{\infty}\geq \frac{\cs\log 
n}{\sqrt n}\Big\}  p_{{}_{X_n(x)}}(u)du\\
 &\qquad \mathop{\longrightarrow}\limits_{\delta\downarrow0}\, 
 \P\Big\{ \|\eps_n(\cdot)-\beta_n(\cdot)X_n(x) \|_{\infty}\geq \frac{\cs\log n}{\sqrt n}\Big\}  p_{{}_{X_n(x)}}(0).
\end{align*} 
Since $|\beta_n(y)|\leq \frac{\cs\log n}{n}$ for $y\in[0,2\pi]$, see part 3 in Lemma \ref{l:aux},
we get
\begin{equation*}
 \P\Big\{ \|\eps_n(\cdot)-\beta_n(\cdot)X_n(x) \|_{\infty}
 \geq \frac{\cs\log n}{\sqrt n}\Big\}  
\leq \P\Big\{ \|\eps_n(\cdot) \|_{\infty}\geq \frac{\cs\log n}{2\sqrt n}\Big\} 
+\P\Big\{ |X_n(x)| \geq \frac{\cs\sqrt{n}}{2}\Big\},
\end{equation*}
which is $o(n)$ since 
the second term is elementary as $X_n(x)\sim N(0,1)$ while the first one follows from part 4 in Lemma \ref{l:aux}. 
The result follows.
\end{proof}
%%%%%%%%%%%%%%%%%%%%%%%%%%%%%%%%%%%%%%%%
\section*{Anciliary computations}
%%%%%%%%%%
%{\bf Lemas viejos de covarianzas:} 
The next Lemma, which proof is omitted, gives the covariances of $X_n$ 
and its derivatives. 
Denote 
\begin{equation} \label{d:rn}
  r_n(x) := \frac1n \sum^n_{k=1} \cos(kx),
\end{equation}
the covariance function of a stationary trigonometric polynomial $T_n$ in \eqref{d:RTP} studied in [GW,AL]. 
\begin{lemma} \label{l:cov}
For $X_n$ and $f_n$ defined as in \eqref{d:Xn} and \eqref{d:fn} respectively, 
we have 
\begin{align*}  
  R_n(x,y) :=  \E(X_n(x)X_n(y)) &= 
   r_n\Big(\frac{\Omega(x)-\Omega(y)}{2}\Big),\\
   R^{(0,1)}_n(x,y) := \E(X_n(x)X'_n(y)) &= 
   -\frac{\omega(y)}{2} r'_n\Big(\frac{\Omega(x)-\Omega(y)}{2}\Big),\\
   R^{(1,1)}_n(x,y) := \E(X'_n(x)X'_n(y)) &= 
   -\frac{\omega(x) \omega(y)}{4} 
   r''_n\Big(\frac{\Omega(x)-\Omega(y)}{2}\Big),\\
   \K_n(x,y) :=\E(f_n(x)f_n(y)) 
  &= r_n\Big(\frac{\Omega(x)-\Omega(y)}{2}\Big) 
  + O\Big(\frac{\log(n)}{n}\Big),\\
  \K^{(0,1)}_n(x,y) := \E(f_n(x)f'_n(y)) 
  &= -\frac{\omega(y)}{2} 
    r'_n\Big(\frac{\Omega(x)-\Omega(y)}{2}\Big)
  + O (1),\\
  \K^{(1,1)}_n(x,y) := \E(f'_n(x)f'_n(y)) 
  &= - \frac{\omega(x)\omega(y)}{4} 
  r''_n\Big(\frac{\Omega(x)-\Omega(y)}{2}\Big) 
  + O (n),\\
   \rho_n(x,y) :=\E(X_n(x)f_n(y)) 
  &= r_n\Big(\frac{\Omega(x)-\Omega(y)}{2}\Big) 
  + O\Big(\frac{\log(n)}{n}\Big),\\
  \rho^{(0,1)}_n(x,y) := \E(X_n(x)f'_n(y)) 
  &= -\frac{\omega(y)}{2} 
    r'_n\Big(\frac{\Omega(x)-\Omega(y)}{2}\Big)
  + O (1),\\
  \rho^{(1,1)}_n(x,y) := \E(X'_n(x)f'_n(y)) 
  &= - \frac{\omega(x)\omega(y)}{4} 
  r''_n\Big(\frac{\Omega(x)-\Omega(y)}{2}\Big) 
  + O (n).
\end{align*}
 where $\Omega$ and $r_n$ are defined 
 in \eqref{e:c-s} and \eqref{d:rn} respectively. 
 In particular, for $x\in[0,2\pi]$ we have
 \[
  \var(F_n(x)) = 1 + O\Big(\frac{\log(n)}{n}\Big), \qquad {\rm and }\
   \var(F_n(x)) = 1 + O\Big(\frac{\log(n)}{n}\Big).
 \]
\end{lemma}

%%%%%%%%%%%%%%%%%%
The next lemma contains some anciliary computations for the proof of Theorem \ref{t:ppal}.
\begin{lemma} \label{l:aux}
With the notation of the proof of Theorem \ref{t:ppal}, we have 
\begin{multicols}{2}
\begin{enumerate}
 \item $|\alpha_n(x)|\leq \cs$ for $n\in\N^*$ and $x\in[0,2\pi]$.
 \item $\Delta=|\det(\var(X_n(x),f_n(x)))|\sim_n \frac{\cs}{n}$.
 \item $|\beta_n(y)|\leq \frac{\cs\log n}{n}$ for $y\in[0,2\pi]$.
 \item $\displaystyle  \P\big\{\|\eps_n\|_{\infty} > \frac{\cs\log n}{\sqrt{n}}\Big\} 
 \leq \frac{1}{n^{\cs\log n}}$.
\end{enumerate}
\end{multicols}
\end{lemma}
\begin{proof}
\begin{enumerate}
\item Its a direct consequence of Lemma \ref{l:cov}.

\item With the notation of Lemma \ref{l:cov}, 
we have 
\begin{align*}
 \Delta &= |\var(X_n(x))\var(f_n(x))- (\E(X_n(x)f_n(x)))^2| \\
 &= |R_n(x,x)\E(\eps_n(x)^2) - (\E(X_n(x)\eps_n(x)))^2|
 \leq \frac{\cs}{n}.
\end{align*}

\item We have
$\displaystyle
 |\beta_n(y)| = |\E(\eps_n(y)X_n(x))| 
 \leq \frac{1}{n}\sum^n_{k=1}O\Big(\frac{1}{k}\Big) 
 \leq\frac{\cs\log n}{n}$.

\item Recall that $\eps_n(x)=\frac{1}{\sqrt n}\sum^n_{k=1} a_kO(\frac{1}{k})+b_kO(\frac{1}{k})$.
Set 
\[
 \sigma^2_n := \sup_x\E|\eps_n(x)|^2 = \E|\eps_n(x_0)|^2
 =\frac{1}{n} \sum^n_{k=1}O\Big(\frac{1}{k}\Big)^2 
 \leq \frac{\cs}{n}.
\]
Besides, 
\begin{multline*}
 \E\Big(\sup_x|\eps_n(x)|\Big) 
 = \frac{1}{\sqrt n} \E\Big[\sup_x \Big|\sum^n_{k=1}a_k O\Big(\frac{1}{k}\Big)
 +b_k O\Big(\frac{1}{k}\Big)\Big|\Big] \\
  \leq \frac{1}{\sqrt n}\sum^n_{k=1}\E|a_k| O\Big(\frac{1}{k}\Big) 
  +\E|b_k| O\Big(\frac{1}{k}\Big) 
 = \frac{1}{\sqrt n}\sum^n_{k=1} O\Big(\frac{1}{k}\Big) 
 \leq\frac{\cs \log n}{\sqrt n}.
\end{multline*}
Now, applying Borell-TIS inequality with $u= \cs\frac{\log n}{\sqrt{n}}-\E\Big(\sup_x|\eps_n(x)|\Big)$, we get
\begin{align*}
 \P\big\{\sup_x|\eps_n(x)| > \frac{\cs\log n}{\sqrt{n}}\Big\} 
 \leq \exp\Big\{ -\frac12 \frac{(\frac{\cs\log n}{\sqrt n})^2}{\frac{\cs}{n}}\Big\} 
 =\frac{1}{n^{\cs\log n}}.
\end{align*}
The result follows.
\end{enumerate}
\end{proof}
%%%
\subsection*{Ackowledgement}
Agradecemos a Juan Pablo Borthagaray y a Guillaume Poly por haber discutido con nosotros sobre algunos aspectos del problema. 
The authors were partially supported by Agencia Nacional de Investigación e Innovación ANII, Uruguay. 
%%%

\end{document}